\documentclass{birkjour}
\usepackage[noadjust]{cite}
\usepackage{xcolor}
\usepackage{tikz}
\RequirePackage[all]{xy}
\usepackage{listings}

\newtheorem{thm}{Theorem}[section]
\newtheorem{cor}[thm]{Corollary}
\newtheorem{lem}[thm]{Lemma}

\theoremstyle{definition}
\newtheorem{defn}[thm]{Definition}
\theoremstyle{remark}
\newtheorem{rem}[thm]{Remark}
\theoremstyle{remark}
\newtheorem{conjecture}[thm]{Conjecture}

\newtheorem{ex}{Example}
\numberwithin{equation}{section}

\newcommand{\BibTeX}{B\kern-0.1emi\kern-0.017emb\kern-0.15em\TeX}
\newcommand{\XYpic}{$\mathrm{X\kern-0.3em\raisebox{-0.18em}{Y}}$-$\mathrm{pic}\,$}

\newcommand{\BZ}{\mathbb{Z}}
\newcommand{\BR}{\mathbb{R}}
\newcommand{\BC}{\mathbb{C}}

\newcommand{\BS}{\mathbb{S}}
\newcommand{\BN}{\mathbb{N}}
\newcommand{\BT}{\mathbb{T}}
\newcommand{\BQ}{\mathbb{Q}}

\newcommand{\ed}{\end{document}}
\begin{document}

%
%
%
%
%
%
%
%

\title[$p$-lck with potential structure on LVMB manifolds]
 {$p$-lck with potential structure on LVMB manifolds}
\author[Bastien Faucard]{Bastien Faucard}
%
\address{%
ENSEA
6 Av. du Ponceau\\
95000 Cergy\\
France\\
ORCID : 0009-0000-5399-4078}
\email{bastien.faucard@ensea.fr}

\subjclass{\\ Primary 32J27, 32L05, 32M99, 32Q15, 32Q60, 32T15, 32V99; \\
Secondary 53A30, 53B35, 53C55, 53C56}
\keywords{LVM, LVMB lck, K\"ahler, lck with potential, $p$-lck structure. }
\date{\today}
\begin{abstract}
In this paper, I present a natural generalization of all the results from \cite{faucard2024lvm} to LVMB manifolds: to summarize, very few LVMB manifolds are lck, and none are lck with potential except for diagonal Hopf manifolds. Moreover, if $N$ is an LVMB manifold with a sufficient number of indispensable coordinates, and under a certain assumption $(H)$ (which may be artificial, as I conjecture) on the localization of the configuration $\Lambda$, there exists a non-compact and non-Kählerian $\BZ^p$-lck with potential cover of $N$ with $p = m-1$. Furthermore, I show that the conjecture stated at the end of \cite{faucard2022lckLVM} (that $p$ is bounded below by $m-1$) is false, by exhibiting examples of LVMB manifolds that are $1$-lck with potential when $m\geq 3$ and $n>2m+1$. This leads to the formulation of a new conjecture: if assumption $(H)$ holds, then $N$ is $1$-lck with potential. Moreover, assumption $(H)$ seems entirely artificial. This conjecture is supported by several examples.
\end{abstract}
\label{page:firstblob}
\maketitle
\newpage
\tableofcontents

\maketitle

\newpage
\section*{Introduction}

The metric aspect of complex manifolds is a fertile field in geometry. For example, Kähler manifolds are complex manifolds that admit a Riemannian metric whose fundamental $2$-form is closed, which grants numerous properties to such objects. A notable example is the deep result known as the Hodge cohomological decomposition theorem. The first historical examples of compact non-Kähler complex manifolds are the Hopf manifolds, introduced in \cite{hopf1930topologie}. Such a manifold is a quotient of $\BC^n\backslash\{0\}$ by a $\BZ$-action generated by a holomorphic contraction. S. López de Medrano generalized this action to an action of $\BC^*\times \BC$ in \cite{lopez1988space} and studied these manifolds in detail in collaboration with A. Verjovsky in \cite{de1997new}. Later, L. Meersseman further extended this to an action of $\BC^*\times \BC^m$, leading to the class of compact complex manifolds known as LVM manifolds \cite{meersseman1998procede}. It was also shown in \cite{meersseman1998procede} that these manifolds, like Hopf manifolds, are non-Kähler.

This naturally raises the question of a more refined metric study of such non-Kähler manifolds. A complex manifold is said to be locally conformally Kähler (we will use the acronym lck) if its metric is locally conformal to a Kähler metric. This class of manifolds was initially introduced by I. Vaisman between the 1970s and 1980s in a series of articles \cite{vaisman1976locally}, \cite{vaisman1980remarkable}, \cite{vaisman1980locally}, \cite{vaisman1982generalized}. It was later studied in detail by L. Ornea and S. Dragomir in 2012 in a book \cite{dragomir2012locally} and in a series of articles by L. Ornea and P. Gauduchon \cite{gauduchon1998locally}, as well as by L. Ornea and M. Verbitsky \cite{ornea2010locally}, \cite{ornea2010topology}, \cite{ornea2014holomorphic}, \cite{ornea2016locally}, \cite{ornea2019hopf}, \cite{ornea2021lee}, \cite{ornea2022non}... Among other contributions to the subject, we may also mention F. Belgun \cite{belgun2000metric}, K. Tsukada \cite{tsukada1999canonical}, and N. Istrati \cite{istrati2019existence}.
One of the rigid properties affecting lck manifolds is that they are not stable under deformations. This result was proved in \cite{ornea2010locally} in 2010. Following this result, the same authors naturally introduced the notion of lck manifolds with potential (where the metric derives from a potential) and showed that such metrics remain stable under deformations.

Hopf manifolds are the archetypical example of lck manifolds with potential. In other hand, in the paper \cite{faucard2024lvm}, the link between lck, lck with potential metrics and LVM manifolds, witch contains Hopf manifolds, is studied. In particular, very few LVM manifolds are lck, except for the diagonal Hopf manifolds.

It is then natural, as I formulated in \cite{faucard2022lckLVM}, to attempt to generalize these results specific to LVM manifolds to larger classes of non-Kähler compact complex manifolds. This article focuses on LVMB manifolds introduced by F. Bosio \cite{bosio2001varietes}, which are a generalization of LVM manifolds. The difference between LVMB and LVM manifolds essentially lies in the fact that Bosio studies, given an action of $\BC^*\times \BC^m$, all possible combinations of subspaces of $\BC^n$ on the complement of which the set of orbits of the action forms a compact complex manifold. The presentation is therefore slightly modified compared to LVM manifolds. Furthermore, the author exactly characterizes LVM manifolds among LVMB manifolds: these are obtained when all the intersections of the interiors of the convex hulls of the possible sub-configurations of the configuration $\Lambda$ are non-empty. This characterization allows for the search for examples that are not LVM manifolds, and thus represent novelties compared to those in \cite{faucard2024lvm}. I would even go further by giving an example borrowed from \cite{cupit2007non} of an LVMB manifold that is not biholomorphic to any LVM manifold. 

Moreover, I present numerous examples of LVMB manifolds to illustrate the results of this paper on the one hand and to demonstrate that the conjecture stated at the end of \cite{faucard2022lckLVM} is false on the other hand. More precisely, there exist LVMB manifolds that are $p$-lck with potential with $p < m-1$. This only makes sense when $m \geq 3$. Going further, the combinatorial study of these examples leads to the formulation of a new conjecture:  

\begin{conjecture}\label{conjecture0}
Let $N$ be an LVMB manifold with $m \geq 1$, $n > 2m+1$, and where the coordinates $1,2,\ldots,m+1$ are essential. Then there exists a non-compact and non-Kählerian $\BZ$-lck cover with potential of $N$.
\end{conjecture}  

In other words, any LVMB manifold satisfying these assumptions that is not a diagonal Hopf manifold (the only LVMB manifolds that are lck with potential) would be $1$-lck with potential.

\vspace{0.5cm}

This article is structured as follows. \\

In the first part, I present the objects involved: lck manifolds, lck manifolds with potential, LVM manifolds, and Hopf manifolds. I also present all the results needed for the rest of the article, citing the contributors.

In the second part, I present the objects that distinguish this article from \cite{faucard2024lvm}: LVMB manifolds, based on the original construction by F. Bosio in \cite{bosio2001varietes}.

In the third part, I propose, in the case where the number of essential coordinates is sufficient, an alternative model for an LVMB manifold, which has the advantage of being a quotient via a discrete action of $\BZ^m$ that replaces the continuous action of $\BC^m$. This part is an adaptation of the corresponding section in \cite{faucard2024lvm} with the same purpose.

The goal of the fourth part is to exactly characterize the diagonal Hopf manifolds among the LVMB manifolds. Again, this is an adaptation of \cite{faucard2024lvm}, however, the proof of Theorem \ref{caraHopf} presents notable differences compared to that of Theorem 3.1 in \cite{faucard2024lvm}. Indeed, two crucial lemmas will be proved beforehand. This theorem essentially states that an LVMB manifold is a diagonal Hopf manifold if and only if $m=1$ and the number of essential coordinates is $2$. 

In the fifth part, I show that under a rationality condition $(K)$ on $\Lambda$, a complex torus of complex dimension $m$ injects naturally into the automorphism group of the LVMB manifold stable under the complex structure. 

In the sixth part, I generalize Theorems 1.9 and 1.11 from \cite{faucard2024lvm} to the case of LVMB manifolds. I obtain the following results (Theorems \ref{thm1} and \ref{thm2}). Let $N$ be an LVMB manifold. If $m=1$ and the number of essential coordinates is $0$ or $1$, then $N$ does not admit an lck metric since it is simply connected. If the number of essential coordinates is $2$, then $N$ is biholomorphic to a diagonal Hopf manifold and thus is lck with potential. If $m \ge 2$, then $N$ is not lck with potential (stability by deformation is crucial for this result). Moreover, if $N$ satisfies the hypothesis $(K)$, then $N$ does not admit an lck metric. Therefore, if there exist LVMB manifolds that are lck, then necessarily $m \ge 2$ and the condition $(K)$ must be violated. Since this condition is dense, such examples, if they exist, would be difficult to uncover.

In the seventh part, I study to what extent an LVMB manifold with $m \ge 2$ and at least $m+1$ essential coordinates (thus admitting a discrete presentation developed in the third part) exhibits an obstruction to admitting an lck metric: the existence of a $\BZ^p$-covering that is lck with potential, non-compact, and non-Kählerian (the key theorem of this part is Theorem \ref{pstructure}). I use, as needed, the terminology of $p$-structure introduced in \cite{faucard2024lvm}.

In the final part, I present several examples of LVMB manifolds that satisfy Theorems \ref{thm1} and \ref{pstructure}: notably, the example from \cite{cupit2007non} where the LVMB manifold with parameters $m=1$ and $n=5$ is not biholomorphic to any LVM manifold and is non-lck according to Theorem \ref{thm1}. I also present examples that refute the conjecture made at the end of \cite{faucard2022lckLVM} regarding the minimality of $p$ in the existence of an lck with potential $p$-structure, which then lead to the formulation of the conjecture \ref{conjecture0}.

Since this search for examples is fundamentally combinatorial in nature, I present in Appendix \ref{Annexe1} all the Python scripts that helped me obtain these examples.

 \section{lck/lck with potential, LVM and Hopf manifolds}

 A locally conformally Kähler (lck) manifold is a complex manifold $N$ whose metric is locally (on an open covering) conformal to a Kähler metric. We can give three equivalent definitions of this notion (see \cite{dragomir2012locally} for the equivalence): the first one reflects the local character, the second provides a necessary and sufficient global condition revealing the closed Lee $1$-form, and finally, the last definition is related to the universal covering.
 
 \vspace{0.1cm}

Let $N = (N, g, J, \Omega)$ be a complex manifold with Riemannian metric $g$, complex structure $J$, and fundamental $2$-form $\Omega$. The manifold $N$ is said to be lck if it satisfies one of the following three equivalent conditions.

\begin{itemize} 
\item \textbf{Local characterization}. There exists an open cover $(U_\alpha)_{\alpha\in A}$ of $N$ and a family of smooth functions $(f\alpha)_{\alpha\in A}$ from $U_\alpha$ to $\BR$ such that, for all $\alpha \in A$, the induced $2$-form 
\begin{equation} \Omega_\alpha = \exp(-f_\alpha)\Omega_{|_{U_\alpha}} \end{equation} 
is closed. 
\item \textbf{Global characterization}. There exists a closed $1$-form $\omega$ such that 
\begin{equation} d\Omega = \omega \wedge \Omega \end{equation} 
This form $\omega$ is called the Lee form. 
\item \textbf{Characterization by covering}. The universal covering $\widetilde{N}$ of $N$ can be endowed with a $2$-form $\widetilde{\Omega}$ such that $(\widetilde{N}, \widetilde{\Omega})$ is a Kähler manifold and the fundamental group $G$ acts on $(\widetilde{N}, \widetilde{\Omega})$ by homotheties (see Theorem 3.31 in \cite{ornea2024principles}). \end{itemize} 
The class of lck metrics is rigid (it is not stable under deformation). It is therefore natural to focus on a certain less rigid subclass of lck manifolds. There are two interesting subclasses: lck manifolds with potential and Vaisman manifolds.

\begin{defn} A lck manifold $N$ is said to be lck with potential if there exists a smooth function $\widetilde{\psi} \colon \widetilde{N} \to \BR_{>0}$ such that $\widetilde{\Omega} = 2i \partial \overline{\partial} \widetilde{\psi} = dd^c \widetilde{\psi}$ and also satisfies that for all $\gamma \in G$, there exists a positive real constant $\chi(\gamma)$ such that $\gamma * \widetilde{\psi} = \chi(\gamma) \widetilde{\psi}$. The function $\widetilde{\psi}$ is called the potential. \end{defn}

\begin{defn} A Vaisman manifold is a lck manifold whose Lee form is parallel: $\nabla \omega = 0$. \end{defn}

Naturally, a Vaisman manifold is a lck manifold with potential, but the converse is generally false.

\vspace{0.1cm}

 It was shown in \cite{ornea2010locally} that an lck metric with potential is stable under deformation, and that a Vaisman metric is not stable under deformation but can deform into an lck metric with potential. We will only use lck manifolds with potential. We are now ready to present the crucial example of lck manifolds with potential (non-Kählerian) announced in the introduction: Hopf manifolds.

 \begin{ex} A Hopf manifold of complex dimension $n \ge 2$ is a quotient of $\BC^n \setminus \{0\}$ by the action of a cyclic group $\langle \alpha \rangle$ (isomorphic to $\BZ$), generated by a holomorphic contraction $\alpha$ ($\alpha(0) = 0$ and for every compact set $K \subset \BC^n$ and every neighborhood $U$ of $0$, there exists a positive integer $q_0$ such that for every integer $q > q_0$, $\alpha^q(K) \subset U$). \ Let $H^n = \BC^n \setminus \{0\} / \langle \alpha \rangle$ be a Hopf manifold. If $\alpha$ is an endomorphism whose eigenvalues have a complex modulus strictly smaller than $1$, we say that $H^n$ is a linear Hopf manifold. If, moreover, $\alpha$ is diagonalizable, we say that $H^n$ is diagonal. In all other cases, we say that $H^n$ is non-linear. It is now known that Hopf manifolds are all compact, complex, non-Kählerian, lck with potential, and the diagonal Hopf manifolds are furthermore Vaisman (see \cite{gauduchon1998locally} for diagonal surfaces, \cite{ornea2016locally} for all linear Hopf manifolds, and \cite{ornea2022non} for all Hopf manifolds). The authors of \cite{ornea2010locally} also demonstrated that for any compact complex manifold $M$ of complex dimension greater than or equal to 3, lck with potential, there exists a holomorphic immersion of $M$ into a linear Hopf manifold. Furthermore, if $M$ is Vaisman, then the linear Hopf manifold is also diagonal. \end{ex}

 \vspace{0.1cm}
 
As a lck manifold with potential (non-Kähler), we have so far only presented Hopf manifolds. Now, a Hopf manifold is generated by a discrete action (an action of $\BZ$). This construction generalizes by considering an action of $\BC^* \times \BC^m$, leading to the manifolds of S. López de Medrano, A. Verjovsky, and L. Meersseman (LVM) \cite{lopez1988space}, \cite{de1997new}, \cite{meersseman1998procede}. We will present the construction of these compact complex (non-Kählerian) manifolds as outlined in \cite{meersseman1998procede}.

\vspace{0.1cm}

Let $m \ge 1$, $n \ge 2m+1$ be two integers and $\Lambda = (\Lambda_1, \ldots, \Lambda_n) \in (\BC^m)^n$ be an $n$-tuple of elements of $\BC^m$ of real rank $2m$ (the maximal rank). We will say that $\Lambda$ is a configuration. For every configuration $\Lambda$, we denote by $\mathcal{H}(\Lambda)$ its real convex hull. Also, let $\Lambda_i = (\lambda_i^1, \ldots, \lambda_i^m) \in \BC^m$. We say that the configuration $\Lambda$ is admissible if it satisfies the following conditions:

\begin{itemize} 
\item the Siegel condition: $0 \in \mathcal{H}(\Lambda)$; 
\item the weak hyperbolicity condition: for all $1 \le i_1 < \cdots < i_{2m} \le n$, $0 \notin \mathcal{H}(\Lambda_{i_1}, \ldots, \Lambda_{i_{2m}})$. 
\end{itemize}
Now, let's define a certain open set of $\BC^n$: 
\begin{align*} 
S_\Lambda = \{z \in \BC^n, 0 \in \mathcal{H}(\Lambda_i)_{i \in I_z}\} 
\end{align*} 
where $I_z = \{i \in \{1, \ldots, n\}, z_i \neq 0\}$ is the set of indices of the non-zero coordinates of $z \in \BC^n$. We have an action of $\BC^* \times \BC^m$ on this open set: 
\begin{align*}
a_\Lambda \colon (\BC^* \times \BC^m) \times S_\Lambda &\rightarrow S_\Lambda \\
(\alpha, T, z) &\mapsto (\alpha z_1 e^{\langle T, \Lambda_1 \rangle}, \ldots, \alpha z_n e^{\langle T, \Lambda_n \rangle}) 
\end{align*} 
where $\langle T, \Lambda_i \rangle = \sum_{j=1}^m T_j \lambda_i^j$.

\begin{thm}[and definition \cite{meersseman1998procede}] 
The quotient $N_\Lambda = S_\Lambda / a_\Lambda$ is called an LVM manifold. It is a compact complex manifold of complex dimension $n - m - 1$. If $n > 2m + 1$, then $N_\Lambda$ is non-symplectic (and hence non-Kählerian). If $n = 2m + 1$, then it is a complex torus $\BT^m$ of complex dimension $m$ (thus a Kähler manifold). \end{thm}

 The following definition will also be necessary: 
 \begin{defn} 
 A coordinate $i \in \{1, \ldots, n\}$ is said to be indispensable if\\
   $0\notin \mathcal{H}(\Lambda_1,\ldots, \Lambda_{i-1},\Lambda_{i+1},\ldots,\Lambda_n)$. Otherwise, we will say that $i$ is eliminable. \end{defn}

Thus, LVM manifolds are non-Kähler (for $n > 2m + 1$), just like Hopf manifolds, from which they are a generalization. In \cite{faucard2024lvm}, we showed that diagonal Hopf manifolds occupy a very special place among LVM manifolds. More specifically, we proved the following theorem:

\begin{thm}\label{caraHopfLVM} 
Let $n > 2m + 1$, $p=n-m-1$, $\Lambda \in (\BC^m)^n$ an admissible configuration with exactly $k$ indispensable coordinates and $N_\Lambda$ the associated LVM manifold. The following assertions are equivalent. 
\begin{enumerate} 
\item[(i)] $N_\Lambda$ is a diagonal Hopf manifold. 
\item[(ii)] $N_\Lambda$ is diffeomorphic to $\BS^{2p-1} \times \BS^1$. 
\item[(iii)] $m = 1$ and $k = 2$. 
\item[(iv)] $N_\Lambda$ has the same homology as $\BS^{2p-1} \times \BS^1$. 
\end{enumerate} \end{thm}

 The next step, in \cite{faucard2024lvm}, was to prove that if an LVM manifold satisfies the following condition $(K_0)$:

\begin{defn}\textbf{(The condition $(K_0)$)}\\ We say that an LVM manifold with an admissible configuration $\Lambda$ satisfies the condition $(K_0)$ if the solution space of the system
\begin{align*} \begin{cases} \sum_{i=1}^n s_i \Lambda_i &= 0 \\ \sum_{i=1}^n s_i &= 0 \end{cases} \end{align*} 
has a rational basis. \end{defn}

Then a complex torus injects into the group of automorphisms stable by the complex structure:

\begin{thm}\label{conditionKToreLVM} Let $N_\Lambda$ be an LVM manifold satisfying the condition $(K_0)$ with $m \ge 1$ and $n \ge 2m+1$. Then there exists a complex torus $\mathbb{T}^m$ in $\mathrm{Aut}(N_\Lambda, J)$. \end{thm}

This theorem allowed us to use an existence criterion for an lck metric given by N. Istrati \cite{istrati2019existence} to prove the following non-existence result.

\begin{thm}\label{ConditionKnonlckLVM} Let $N_\Lambda$ be an LVM manifold with $n > 2m + 1$. Suppose that $N_\Lambda$ is not biholomorphic to a diagonal Hopf manifold. If $m = 1$ or $m \ge 2$ and $N_\Lambda$ satisfies the condition $(K_0)$, then there does not exist an lck metric on $N_\Lambda$. \end{thm}

\begin{rem} Since lck metrics are not stable under deformation, it is possible that there exist LVM manifolds admitting an lck metric (other than diagonal Hopf manifolds). However, to obtain such examples, it is necessary for $m$ to be greater than $2$ and for $\Lambda$ to not satisfy the condition $(K_0)$, which is a dense condition: this leaves relatively little room. We do not have such an example to date, but to complete this remark, let us add the following theorem, which follows from the stability of the lck with potential metric structure under deformations: there does not exist any (except for diagonal Hopf manifolds) LVM manifold that is lck with potential. \end{rem}

\begin{thm}\label{CorollairecoolLVM} There exists an lck with potential metric on an LVM manifold if and only if it is biholomorphic to a diagonal Hopf manifold. \end{thm}

In \cite{faucard2024lvm}, following this result of the non-existence of an lck with potential metric, I study the existence of $\BZ^p$-coverings that are lck with potential. More specifically, I prove the following theorem:

\begin{thm} Let $N$ be an LVM manifold with an admissible configuration $\Lambda$ whose first $m+1$ coordinates are indispensable, satisfying the following condition:\ There exists a basis $\mathcal{B}=(f_1, \ldots, f_m)$ of $\BZ^m$ and $j \in {1, \ldots, m}$ such that \begin{align*} |\Gamma(f_j)_r|<1, \forall r\in {m+2, \ldots, n}. \end{align*} See \cite{faucard2024lvm} for the construction of the $\Gamma(f_j)_r$ or equation \ref{Gamma} from section \ref{contdiscret} of this article. Then there exists a non-compact, non-Kählerian $\BZ^{m-1}$-covering lck with potential of $N$. Moreover, this result remains stable under deformation (see Remark 6.3 of \cite{faucard2024lvm} on this topic). \end{thm}

In the next section, we will present a generalization of the LVM manifolds due to F. Bosio \cite{bosio2001varietes}: the LVMB manifolds. These include the LVM manifolds and are also compact, complex, non-Kählerian manifolds (if $n > 2m+1$).

\section{LVMB manifolds}

Let $n \geq 1$ be an integer and $m \geq 1$ another integer. A fundamental set is a non-empty set of subsets of $\{1, \ldots, n\}$, where all elements of the set have cardinality $2m+1$. This definition is relative to $m$, so if $n < 2m+1$, then every fundamental set is empty. For the rest of the discussion, we assume that $n \geq 2m+1$. Let $\epsilon$ be a fundamental set. An element $\sigma \in \epsilon$ is called a fundamental part. A subset $P$ of $\{1, \ldots, n\}$ is said to be acceptable if it contains a fundamental part $\sigma \in \epsilon$. Let $\mathcal{A}$ denote the subset of $\mathcal{P}(\{1, \ldots, n\})$ consisting of acceptable subsets. A coordinate $i$ is called indispensable if it is contained in all the fundamental parts. 

\vspace{0.1cm}

Let $\epsilon$ be a fundamental set. Consider an $n$-tuple $\Lambda \in (\BC^m)^n$ of vectors in $\BC^m$, which is called a configuration. The pair $(\epsilon, \Lambda)$ will be called a studyable system if, for every fundamental subset $\sigma \in \epsilon$, $(\Lambda_i)_{i \in \sigma}$ forms a real affine frame of $\BC^m = \BR^{2m}$.

Let $(\epsilon, \Lambda)$ be a studyable system. Define
\begin{align*}
S^\epsilon=\{z\in \BC^n,I_z\in \mathcal{A}\},\end{align*}
and consider the action
\begin{align*}
a_\Lambda\colon (\BC^*\times \BC^m)\times S^\epsilon&\rightarrow S^\epsilon\\
(\alpha,T,z)&\mapsto (\alpha z_1 e^{<\Lambda_1,T>},\ldots,\alpha z_n e^{<\Lambda_n,T>}).
\end{align*}
Let $N^\epsilon_\Lambda$ denote the quotient of $S^\epsilon$ by the action $a_\Lambda$. Bosio showed in \cite{bosio2001varietes} that $N^\epsilon_\Lambda$ is a compact complex manifold if and only if the studyable system $(\epsilon, \Lambda)$ satisfies the imbrication condition and the $(PEUR)$: principle of existence and uniqueness of the replacer.

\begin{defn} We say that a fundamental set $\epsilon$ satisfies the \ $(PER)$: principle of existence of the replacer if, for every fundamental part $\sigma \in \epsilon$ and every $k \in \{1, \ldots, n\}$, there exists $k' \in \sigma$ such that $(\sigma \setminus \{k'\}) \cup \{k\} \in \epsilon$. We say that a fundamental set satisfying the $(PER)$ satisfies the $(PEUR)$ if the $k'$ is unique. \end{defn}

\begin{defn} We say that an studyable system $(\epsilon, \Lambda)$ satisfies the imbrication condition if for every pair $(\sigma, \tau)$ of fundamental parts, the real convex hulls $\mathcal{H}(\Lambda_i)_{i \in \sigma}$ and $\mathcal{H}(\Lambda_j)_{j \in \tau}$ have non-disjoint interiors. \end{defn}

A studyable system that satisfies the imbrication condition and the $(PEUR)$ is called a good system, and in this case, the compact complex manifold $N_\Lambda^\epsilon$ is called an LMVB manifold.\\

Bosio, by generalizing Meersseman’s proof, showed that if $n > 2m + 1$, the LVMB manifolds are non-symplectic \cite{bosio2001varietes}. Furthermore, he showed that any LVM manifold can be realized as an LVMB manifold, but the converse is not true in general. More precisely:

\begin{thm}\label{caractLVM} Let $\Lambda \in (\BC^m)^n$ be an admissible configuration. Let $\theta$ be the set of points in $\BC^m$ that do not belong to any part of $\Lambda$ with $2m$ elements. Then, $(\epsilon, \Lambda)$ is the good system of an LVM manifold if and only if there exists a bounded connected component $O \subset \theta$ such that 
\begin{align*} 
\epsilon = \{\sigma \in \mathcal{P}(\{1, \ldots, n\}), |\sigma| = 2m + 1, O \subset \mathcal{H}(\Lambda_i)_{i \in \sigma}\}. 
\end{align*} 
\end{thm}

\section{From continous to discret}\label{contdiscret}

As with LVM manifolds, an LVMB manifold can be presented using a discrete action of $\BZ^m$ (if it has sufficiently many indispensable coordinates).

\begin{thm}\label{continu_discret}
Let $N^\epsilon_\Lambda$ be an LVMB manifold with $m \ge 1$, $n > 2m + 1$, and $k \ge m + 1$ indispensable coordinates, for which we can assume that the first $(m + 1)$ coordinates, $1, 2, \ldots, m + 1$, are included. Then there exists an open set $C \subset \BC^{n-m-1}$ and an action $d_\Lambda$ of $\BZ^m$ on $C$ such that
\begin{align*}
N^\epsilon_\Lambda=C/d_\Lambda.\end{align*}
\end{thm}
 
\begin{proof}
First, according to the hypotheses on $\epsilon$, the open set $S^\epsilon$ is written as
\begin{align*}
S^\epsilon = (\BC^*)^{m+1} \times C,
\end{align*}
where $C$ is an open set in $\BC^{n-m-1}$. The reasoning that follows is exactly the same as in the "from continuous to discrete" section of \cite{faucard2024lvm}. The only new point is the proof that the matrix
\begin{align*}
A = \begin{pmatrix} \Lambda_2 - \Lambda_1 \\ \vdots \\ \Lambda_{m+1} - \Lambda_1 \end{pmatrix}
\end{align*}
is invertible. This naturally arises in the reasoning for constructing the discrete action of $\BZ^m$ in the case of LVM manifolds. In this case, the weak hyperbolicity condition was used. Here, we use the fact that the system $(\epsilon, \Lambda)$ is studyable and that $1, 2, \ldots, m+1$ are indispensable. Indeed, by assumption, we have an $\BR$-basis of $\BC^m$: $(\Lambda_j - \Lambda_1)_{j \in \sigma \setminus \{1\}}$ for all $\sigma \in \epsilon$. Thus, since $2, 3, \ldots, m+1 \in \sigma$, the family $(\Lambda_j - \Lambda_1)_{2 \le j \le m+1}$ is $\BR$-free because it is a subfamily of a $\BR$-basis. Therefore, the matrix $A$ is indeed invertible. To have the necessary notations, we will explicitly write out the action of $\BZ^m$:
\begin{align*}
d_\Lambda \colon \BZ^m \times C &\rightarrow C \\
(l, w_{m+2}, \ldots, w_n) &\mapsto \mathrm{Diag}(\Gamma(l)_{m+2}, \ldots, \Gamma(l)_n) \cdot \begin{pmatrix} w_{m+2} \\ \vdots \\ w_n \end{pmatrix}
\end{align*}
where, for each $r \in \{m+2, \ldots, n\}$,
\begin{equation}\label{Gamma}
\Gamma(l)_r = e^{\langle \Lambda_r - \Lambda_1, 2i \pi A^{-1} \cdot l \rangle}.
\end{equation}
In \cite{faucard2024lvm}, it is shown (just before Theorem 2.1) that $N^\epsilon_\Lambda$ is biholomorphic to $C/d_\Lambda$. If $\mathcal{B} = (f_1, \ldots, f_m)$ is a basis of $\BZ^m$, the action of $\BZ^m$ defining $N^\epsilon_\Lambda$ is generated by
\begin{align*}
\alpha_1^\mathcal{B}, \ldots, \alpha_m^\mathcal{B},
\end{align*}
where for each $j \in \{1, \ldots, m\}$,
\begin{equation}\label{alpha}
\alpha_j^\mathcal{B} = \mathrm{Diag}(\Gamma(f_j)_{m+2}, \ldots, \Gamma(f_j)_n).
\end{equation}
\end{proof}

\section{Diagonal Hopf manifolds ans LVMB manifolds}

This section is devoted to the proof of the following theorem, a generalization of Theorem 3.1 in \cite{faucard2024lvm}.

\begin{thm}\label{caraHopf}
Let $n>2m+1$ and $N^\epsilon_\Lambda$ be an LVMB manifold of type $(n,m,k)$. Denote by $p=n-m-1$. The following three assertions are equivalent.
\begin{enumerate}
\item $N^\epsilon_\Lambda$ is a diagonal Hopf manifold.
\item $N^\epsilon_\Lambda$ is diffeomorphic to $\BS^{2p-1}\times \BS^1$.
\item $m=1$ and $k=2$.
\end{enumerate}
\end{thm}
The proof of this theorem relies on the following two lemmas.

\begin{lem}\label{lem1}
Let $\epsilon$ be a fundamental set with $m=1$ where the coordinates $1$ and $2$ are indispensable. Then
\begin{align*}
\epsilon = \{(12j), j = 3, \ldots, n\}.
\end{align*}
\end{lem}
\begin{proof}
We will prove the equality of sets by double inclusion. Let $\sigma \in \epsilon$. Since $1$ and $2$ are indispensable, we have $1,2 \in \sigma$, and since $m=1$, $|\sigma|=3$, so there exists $j \in \{3, \ldots, n\}$ such that $\sigma = (12j)$. Thus, we have the direct inclusion. Moreover, from this fact, we know that $\epsilon = \{(12j), j \in E\}$ where $E$ is a subset of $\{3, \ldots, n\}$. Let $E = \{j_1, \ldots, j_l\}$ with $j_i \in \{3, \ldots, n\}$, and thus $l \le n-2$. Let $j \in \{3, \ldots, n\}$. We will show that $\sigma = (12j) \in \epsilon$ for the reverse inclusion. If $j \in E$, then $\sigma = (12j) \in \epsilon$. Suppose that $j \notin E$. Let $j_i \in E$, and let $\tau = (12j_i) \in \epsilon$. By the $(PEUR)$ property, there exists a unique $p \in \tau$ such that 
\begin{align*}
(\tau \backslash \{p\}) \cup \{j\} \in \epsilon.
\end{align*}
Since $1,2$ are indispensable, we have $p = j_i$ and $\sigma = \tau \backslash \{j_i\} \cup \{j\} \in \epsilon$. The reverse inclusion is thus established. In particular, $E = \{3, \ldots, n\}$.
\end{proof}

\begin{lem}\label{lem2}
Let $(\epsilon, \Lambda)$ be a good system with $m = 1$ where the coordinates $1$ and $2$ are indispensable. Let $D$ be the line in the complex plane $\BC$ passing through $\Lambda_1$ and $\Lambda_2$. Then the sign of $\mathrm{Im}\Big(\frac{\Lambda_1 - \Lambda_j}{\Lambda_2 - \Lambda_1}\Big)$ is constant for $j \in \{3, \ldots, n\}$ if and only if all the $\Lambda_j$, $j \in \{3, \ldots, n\}$, are on the same side of $D$.
\end{lem}

\begin{proof}
First, observe that $\Lambda_2 \neq \Lambda_1$. Indeed, since $(\epsilon, \Lambda)$ is a good system, for every $\sigma \in \epsilon$, $(\Lambda_i)_{i \in \sigma}$ is a real affine frame of $\BC$, and since $1, 2 \in \sigma$, we must have $\Lambda_1 \neq \Lambda_2$.\\
Let $\mathrm{sign}(x)$ denote the function that equals $1$ if the real number $x$ is positive and $-1$ if it is negative (this function is not defined at $0$, but we do not need this case here).\\
Let for each $j \in \{3, \ldots, n\}$, $\Lambda_j = a_j + i b_j$ with $a_j, b_j \in \BR$. Suppose first that $a_1 = a_2$ (i.e., that the line $D$ is vertical). We have
\begin{align*}
\mathrm{sign}\Big(\mathrm{Im}\Big(\frac{\Lambda_1 - \Lambda_j}{\Lambda_2 - \Lambda_1}\Big)\Big)
&= \mathrm{sign}\Big(\mathrm{Im}\Big(\frac{\Lambda_1 - \Lambda_j}{i(b_2 - b_1)}\Big)\Big) \\
&= \mathrm{sign}\Big(\mathrm{Im}\Big(\frac{i(\Lambda_1 - \Lambda_j)}{b_1 - b_2}\Big)\Big) \\
&= \mathrm{sign}\Big(\frac{a_1 - a_j}{b_1 - b_2}\Big) \\
&= \mathrm{sign}(b_1 - b_2)\mathrm{sign}(a_1 - a_j).
\end{align*}
This sign is constant in $j$ if and only if the sign of $a_1 - a_j$ is constant in $j$. In other words, this is equivalent to the fact that all the $a_j$ are strictly less than (or strictly greater than) $a_1$, meaning that all the $\Lambda_j$ are on the same side of $D$. Now, suppose that $a_1 \neq a_2$. The equation of the line $D$ is $y = \alpha x + \beta$ with $\alpha = \frac{b_2 - b_1}{a_2 - a_1}$ and $\beta = b_1 - \alpha a_1$. We have
\begin{align*}
\mathrm{sign}\Big(\mathrm{Im}\Big(\frac{\Lambda_1 - \Lambda_j}{\Lambda_2 - \Lambda_1}\Big)\Big)
&= \mathrm{sign}\Big(\mathrm{Im}\Big(\frac{a_1 - a_j + i(b_1 - b_j)}{a_2 - a_1 + i(b_2 - b_1)}\Big)\Big) \\
&= \mathrm{sign}\Big(\frac{-(a_1 - a_j)(b_2 - b_1) + (a_2 - a_1)(b_1 - b_j)}{|\Lambda_2 - \Lambda_1|^2}\Big) \\
&= \mathrm{sign}\Big(-(a_1 - a_j)(b_2 - b_1) + (a_2 - a_1)(b_1 - b_j)\Big).
\end{align*}
We can then force a factorization by $(a_2 - a_1)$:
\begin{align*}
\mathrm{sign}\Big(\mathrm{Im}\Big(\frac{\Lambda_1 - \Lambda_j}{\Lambda_2 - \Lambda_1}\Big)\Big)
&= \mathrm{sign}\Big((a_2 - a_1)\Big(\frac{-(a_1 - a_j)(b_2 - b_1)}{a_2 - a_1} + b_1 - b_j\Big)\Big) \\
&= \mathrm{sign}(a_2 - a_1)\mathrm{sign}\Big(-\alpha(a_1 - a_j) + b_1 - b_j\Big) \\
&= \mathrm{sign}(a_2 - a_1)\mathrm{sign}(\beta + \alpha a_j - b_j).
\end{align*}
If all the $\Lambda_j$ are on the same side of $D$, the sign above is constant, and conversely.
\end{proof}

We can now begin the proof of theorem \ref{caraHopf}.
\begin{proof}
\begin{enumerate}
\item \textit{First step.} Suppose that $N^\epsilon_\Lambda$ is a diagonal Hopf manifold, then by Theorem 3.1 \cite{faucard2024lvm} or Theorem 5.3.1 from \cite{faucard2022lckLVM}, we directly get $(2)$ and $(3)$, as in this case $N^\epsilon_\Lambda$ is an LVM manifold. \\
In the next two steps, we assume $(3)$: $k=2$ and $m=1$, with coordinates $1$ and $2$ being indispensable. 
\item \textit{Second step.} Recall that
\begin{align*}
S^\epsilon = \{z \in \BC^n, I_z \in \mathcal{A}\},
\end{align*}
and show that 
\begin{align*}
S^\epsilon = (\BC^*)^2 \times \BC^{n-2} \backslash \{0\}.
\end{align*}
We do this by double inclusion. Let $z \in (\BC^*)^2 \times \BC^{n-2} \backslash \{0\}$, we need to show that $z \in S^\epsilon$, that is, that $I_z$ is acceptable, hence contains a fundamental part. By assumption, there exists $j \in \{3, \ldots, n\}$ such that $(12j) \in I_z$. By Lemma \ref{lem1}, $(12j) \in \epsilon$: hence $z \in S^\epsilon$. \\
Conversely, if $z \in S^\epsilon$, then $I_z$ is an acceptable part, so it contains a fundamental part (of the form $\sigma = (12j)$ with $j \in \{3, \ldots, n\}$ by Lemma \ref{lem1}). This shows that $z_1 \neq 0$, $z_2 \neq 0$, and there exists $j \in \{3, \ldots, n\}$ such that $z_j \neq 0$, thus $z \in (\BC^*)^2 \times \BC^{n-2} \backslash \{0\}$. We have the equality 
\begin{align*}
S^\epsilon = (\BC^*)^2 \times \BC^{n-2} \backslash \{0\}.
\end{align*}
\item \textit{Third step.} Apply Theorem \ref{continu_discret} with $m = 1$ and $k = 2$ (where we can assume that coordinates $1$ and $2$ are indispensable). The second step has shown that $C = \BC^{n-2} \backslash \{0\}$. The action $d_\Lambda$ is an action of $\BZ$. Let $l \in \BZ^*$ be a non-zero integer (thus a generator of $\BZ$), then the action $d_\Lambda$ is generated by $\alpha_1^l = \mathrm{Diag}(\Gamma(l)_3, \ldots, \Gamma(l)_n)$ with 
\begin{align*}
\Gamma(l)_r = e^{<\Lambda_r - \Lambda_1, 2i\pi A^{-1} \cdot l>} = e^{(\Lambda_r - \Lambda_1) \times \frac{2i\pi l}{\Lambda_2 - \Lambda_1}}.
\end{align*}
We need to show that $\alpha_1^l$ is a holomorphic contraction for some $l \in \BZ^*$. In other words, we need to show that there exists $l \in \BZ^*$ such that 
\begin{align*}
|\Gamma(l)_r| < 1
\end{align*}
for all $r \in \{3, \ldots, n\}$. This translates to
\begin{align*}
\Big| e^{2i\pi l \frac{\Lambda_1 - \Lambda_r}{\Lambda_2 - \Lambda_1}} \Big| < 1 &\Leftrightarrow e^{-2\pi l \mathrm{Im}\Big(\frac{\Lambda_1 - \Lambda_r}{\Lambda_2 - \Lambda_1}\Big)} < 1 \\
&\Leftrightarrow -2\pi l \mathrm{Im} \Big(\frac{\Lambda_1 - \Lambda_r}{\Lambda_2 - \Lambda_1}\Big) < 0 \\
&\Leftrightarrow l \mathrm{Im} \Big(\frac{\Lambda_1 - \Lambda_r}{\Lambda_2 - \Lambda_1}\Big) > 0.
\end{align*}
Since $l$ can be chosen to be positive or negative, it is sufficient to show that $\mathrm{Im} \Big(\frac{\Lambda_1 - \Lambda_r}{\Lambda_2 - \Lambda_1}\Big)$ has the same sign for all $r \in \{3, \ldots, n\}$. By Lemma \ref{lem2}, this condition is equivalent to all the $\Lambda_r$ (for $r \in \{3, \ldots, n\}$) being on the same side of the line $D$ passing through $\Lambda_1$ and $\Lambda_2$. Suppose by contradiction that there exist distinct $i_1, i_2 \in \{3, \ldots, n\}$ such that $\Lambda_{i_1}$ and $\Lambda_{i_2}$ are not on the same side of $D$. By Lemma \ref{lem1}, the parts $(12i_1)$ and $(12i_2)$ are fundamental. This situation violates the imbrication condition since the interiors of the convex hulls $\mathcal{H}(\Lambda_1, \Lambda_2, \Lambda_{i_1})$ and $\mathcal{H}(\Lambda_1, \Lambda_2, \Lambda_{i_2})$ are disjoint (see Figure \ref{fig1}). Therefore, all the $\Lambda_r$ are on the same side of $D$: thus $(1)$.
\end{enumerate}

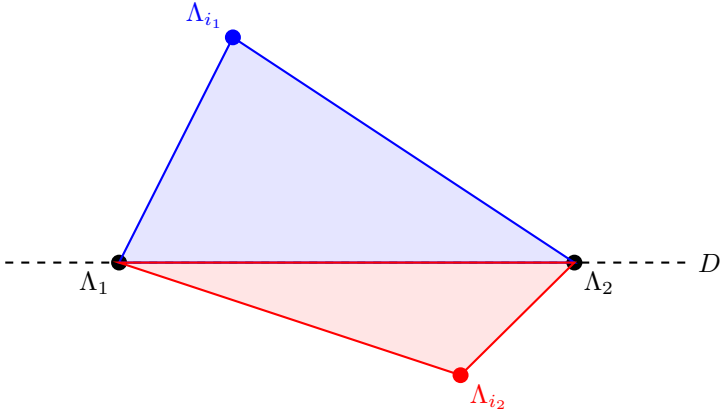
\begin{figure}[h!]

\begin{tikzpicture}[scale=1.5]

\coordinate (Lambda1) at (0, 0);
\coordinate (Lambda2) at (4, 0);
\coordinate (Lambda_i1) at (1, 2);
\coordinate (Lambda_i2) at (3, -1);

\draw[thick, dashed, black] (-1, 0) -- (5, 0) node[right] {$D$};

\fill[blue!20, opacity=0.5] (Lambda1) -- (Lambda2) -- (Lambda_i1) -- cycle;
\fill[red!20, opacity=0.5] (Lambda1) -- (Lambda2) -- (Lambda_i2) -- cycle;

\fill[black] (Lambda1) circle (2pt) node[below left] {$\Lambda_1$};
\fill[black] (Lambda2) circle (2pt) node[below right] {$\Lambda_2$};
\fill[blue] (Lambda_i1) circle (2pt) node[above left] {$\Lambda_{i_1}$};
\fill[red] (Lambda_i2) circle (2pt) node[below right] {$\Lambda_{i_2}$};

\draw[thick, blue] (Lambda1) -- (Lambda2) -- (Lambda_i1) -- cycle;
\draw[thick, red] (Lambda1) -- (Lambda2) -- (Lambda_i2) -- cycle;


\end{tikzpicture}
   \caption{$\Lambda_{i_1}, \Lambda_{i_2}$ and the line passing through $\Lambda_1$ and $\Lambda_2$.}
   \label{fig1}
   \end{figure}
\end{proof}

\section{Existence of a complex torus in $\mathrm{Aut}_J(N^\epsilon_\Lambda)$}

\begin{defn}
Let $(\epsilon,\Lambda)$ be a good system. We say that it (or the associated LVMB manifold $N_\Lambda^\epsilon$) satisfies the condition $(K)$ (we will also write $\Lambda \in (K)$ or $N_\Lambda^\epsilon \in (K)$) if there exists a $\BR$-automorphism $\Phi$ of $\BC^m$ such that for every $j \in \{1, \ldots, n\}$, $\Phi(\Lambda_j) \in \BZ[i]^m$.
\end{defn}

Let us make two remarks.

\begin{rem}
Recall that the condition $(K_0)$ is the existence of a rational basis for the system
\begin{align*}
(S) \colon \begin{cases}
\sum_{j=1}^n s_j \Lambda_j &= 0 \\
\sum_{j=1}^n s_j &= 0
\end{cases}.
\end{align*}
The condition $(K)$ implies the condition $(K_0)$. Indeed, let $r$ be the rank of the system $(S)$ viewed with real coefficients. This rank is always equal to $2m$ since the system $(\epsilon,\Lambda)$ is studyable. If the condition $(K)$ is satisfied, there exists an $\BR$-automorphism $\Phi$ of $\BC^m$ such that $\Phi(\Lambda_j) \in \BZ[i]^m$ for every $j$. Applying $\Phi$ to the previous system, it is then equivalent to
\begin{align*}
(S_\Phi) \colon \begin{cases}
\sum_{j=1}^n s_j \Phi(\Lambda_j) &= 0 \\
\sum_{j=1}^n s_j &= 0
\end{cases}.
\end{align*}
Moreover, $(S_\Phi)$ is a system with coefficients in $\BZ \subset \mathbb{Q}$, its rank over $\mathbb{Q}$ is equal to $r$ since the rank is the largest non-zero minor and the determinant is independent of the base field. It follows that the system $(S_\Phi)$ satisfies $n > r = 2m$, so it admits a rational solution basis, and the same holds for the equivalent system $(S)$. Therefore, the condition $(K_0)$ is satisfied.
\end{rem}

\begin{rem}
Suppose that $\Lambda \in (\mathbb{Q}[i]^m)^n$, then the condition $(K)$ is satisfied (it suffices to define the automorphism by multiplying by the product of all the denominators of the real and imaginary parts of the coordinates of the $\Lambda_j$) and thus the condition $(K_0)$ is also satisfied. Since $\mathbb{Q}$ is dense in $\BR$, these two conditions are dense: more precisely, if $\epsilon$ is a fixed fundamental set, the set
\begin{align*}
\{(\epsilon,\Lambda), \Lambda \in (K)\}
\end{align*}
is dense in the set of good systems with fundamental set $\epsilon$.
\end{rem}

 We will now prove the following theorem, which is an adaptation of Theorem 4.2 from \cite{faucard2024lvm} to LVMB manifolds.

\begin{thm}\label{ToreLVMB}
Let $ N^\epsilon_\Lambda $ be an LVMB manifold satisfying the condition $(K)$ with $ n \ge 2m + 1 $ and $ m \geq 1 $. Then there exists a complex torus $ \BT^m $ of complex dimension $ m $ in the automorphism group of $ N^\epsilon_\Lambda $ that is preserved by the complex structure $ J $.
\end{thm}

\begin{proof}
Consider the action
\begin{align*}
\eta\colon \BC^*\times \BC^{2m}\times S^\epsilon &\rightarrow S^\epsilon \\
((\alpha,T,U),z) &\mapsto \Big(\alpha z_j e^{<\Lambda_j,T>+<\mathrm{Re}(\Lambda_j),U>}, j=1,\ldots,n\Big),
\end{align*}
and its group of ineffectivity:
\begin{align*}
\Gamma = \{(\alpha,T,U) \in \BC^*\times \BC^{2m}, \forall z \in S^\epsilon, \eta(\alpha,T,U) = z\}.
\end{align*}
This action naturally induces another one of
\begin{align*}
\BC^m/\Gamma \approx \Big(((\BC^*\times \BC^{2m})/\Gamma)/(\BC^*\times \BC^m)\Big)
\end{align*}
on $N_\Lambda^\epsilon = S^{\epsilon}/a_\Lambda$ (the symbol $\approx$ here means that both parts are biholomorphic). Now, if $(\alpha,T,U) \in \Gamma$ and $z \in (\BC^*)^n$, we have
\begin{align*}
\alpha e^{<\Lambda_j,T>+<\mathrm{Re}(\Lambda_j),U>} = 1.
\end{align*}
Next, for each equation, divided by the one where $j=1$ and taken logarithmically, we get
\begin{align*}
<\Lambda_j-\Lambda_1,T>+<\mathrm{Re}(\Lambda_j)-\mathrm{Re}(\Lambda_1),U> = 2i\pi k_j,
\end{align*}
for some $(k_2,\ldots,k_n) \in \BZ^{n-1}$. Let $P = T + U$ and $Q = iT$, this yields
\begin{align*}
<\mathrm{Re}(\Lambda_j-\mathrm{Re}(\Lambda_1),P> + <\mathrm{Im}(\Lambda_j)-\mathrm{Im}(\Lambda_1),Q> = 2i\pi k_j.
\end{align*}
Finally, separate the real and imaginary parts of this system:
\begin{align*}
<\mathrm{Re}(\Lambda_j-\mathrm{Re}(\Lambda_1),\mathrm{Re}(P)> + <\mathrm{Im}(\Lambda_j)-\mathrm{Im}(\Lambda_1),\mathrm{Re}(Q)> &= 0, \\
<\mathrm{Re}(\Lambda_j-\mathrm{Re}(\Lambda_1),\mathrm{Im}(P)> + <\mathrm{Im}(\Lambda_j)-\mathrm{Im}(\Lambda_1),\mathrm{Im}(Q)> &= 2\pi k_j.
\end{align*}
The first system (where we consider the real part) is a system of $n-1$ equations with $2m$ variables. Let $\sigma \in \epsilon$, since $(\epsilon,\Lambda)$ is studyable, $(\mathrm{Re}(\Lambda_j))_{j \in \sigma}$ is an affine real basis of $\BR^m$. If $1 \in \sigma$, then $(\mathrm{Re}(\Lambda_j-\Lambda_1))_{j \in \sigma \backslash \{1\}}$ is an $\BR$-basis of $\BR^m$. If $1 \notin \sigma$, then by $(PEUR)$, there exists a unique replacement $k \in \sigma$ such that $\tau = (\sigma \backslash \{k\}) \cup \{1\} \in \epsilon$. Then, $(\mathrm{Re}(\Lambda_j-\Lambda_1))_{j \in \tau \backslash \{1\}}$ is an $\BR$-basis of $\BR^m$. Proceeding similarly for the imaginary parts, we find that there exists $\sigma \in \epsilon$ containing $1$ such that
\begin{align*}
(\Lambda_j-\Lambda_1)_{j \in \sigma \backslash \{1\}} = (\mathrm{Re}(\Lambda_j-\Lambda_1), \mathrm{Im}(\Lambda_j-\Lambda_1))_{j \in \sigma \backslash \{1\}}
\end{align*}
is an $\BR$-basis of $\BC^m$. This demonstrates that there exists a submatrix of this system with maximal real rank $2m$: we then have
\begin{align*}
\mathrm{Re}(P) = \mathrm{Re}(Q) = 0.
\end{align*}
Now, consider the second part (where we only use the imaginary part):
\begin{align*}
I_j = <\mathrm{Re}(\Lambda_j-\Lambda_1), \mathrm{Im}(P)> + <\mathrm{Im}(\Lambda_j-\Lambda_1), \mathrm{Im}(Q)> = 2\pi k_j.
\end{align*}
By the same reasoning as before, there exists $\sigma \in \epsilon$ containing $1$ such that $(\Lambda_j-\Lambda_1)_{j \in \sigma \backslash \{1\}}$ is an $\BR$-basis of $\BC^m$. Thus, for each $(k_j)_{j \in \sigma \backslash \{1\}}$, the $(I_j)_{j \in \sigma \backslash \{1\}}$ are fixed. The remaining equations are those with $j \in \{2,\ldots,n\} \backslash \sigma$. This is where the condition $(K)$ comes into play. Since $(K)$ implies $(K_0)$, there exists a combination
\begin{align*}
\Lambda_j - \Lambda_1 = \sum_{l \in \sigma \backslash \{1\}} a_{j,l} (\Lambda_l - \Lambda_1),
\end{align*}
where the $a_{j,l}$ can be assumed to be in $\BZ$ (by multiplying these rationals by a common integer). And then, for these same $j$:
\begin{align*}
I_j = 2\pi k_j,
\end{align*}
where
\begin{align*}
k_j = \sum_{l \in \sigma \backslash \{1\}} a_{j,l} k_l \in \BZ
\end{align*}
is determined. To summarize, the stabilizer under the action of $\BC/\Gamma$ at a point $z$ is identified with $\BZ^{2m}$. Moreover, the elements of this stabilizer fix each point of $S^\epsilon$. This means that $\Gamma$ is identified with $\BZ^{2m}$. Thus:
\begin{align*}
\BC^m/\Gamma = \BC^m/\BZ^{2m} = \BT^m.
\end{align*}
Furthermore, the action of this complex torus of complex dimension $m$, being totally effective by construction, induces the injection of $\BT^m$ into the automorphism group of $N^\epsilon_\Lambda$ that is preserved by $J$: $\mathrm{Aut}(N^\epsilon_\Lambda, J)$.
\end{proof}

\section{lck metrics and LVMB manifolds}

The conclusions are then, by what has been previously shown, similar to those in \cite{faucard2024lvm}. More precisely, Theorems \ref{ConditionKnonlckLVM} and \ref{CorollairecoolLVM} naturally generalize to LVMB manifolds. Here they are in the case of LVMB manifolds.

\begin{thm}\label{thm1}
Let $N^\epsilon_\Lambda$ be an LVMB manifold of type $(n,m,k)$ with $n>2m+1$. Suppose that $N^\epsilon_\Lambda$ is not biholomorphic to a diagonal Hopf manifold (in other words, according to Theorem \ref{caraHopf}, $(m,k)\neq (1,2)$). If $m=1$, or if $m>1$ and $N^\epsilon_\Lambda$ satisfies condition $(K)$, then there does not exist lck metric on $N^\epsilon_\Lambda$.
\end{thm}
\begin{proof}
Let $N^\epsilon_\Lambda$ be an LVMB manifold that is not biholomorphic to a diagonal Hopf manifold. First, suppose that $m=1$. The number of indispensable points is less than $2m+1=3$. Suppose that $k=3$, then all the fundamental parts are identical. The set $\epsilon$ satisfying this condition does not satisfy the $(PEUR)$. Therefore, $k\neq 3$. According to Theorem \ref{caraHopf}, we cannot have $k=2$. The remaining cases are $k=0$ and $k=1$. But for both of these cases, $N^\epsilon_\Lambda$ becomes simply connected, hence non-lck.\\
Now suppose that $m\ge 2$ and that $N^\epsilon_\Lambda$ satisfies condition $(K)$. By Theorem \ref{ToreLVMB}, there exists a complex torus $\BT^m$ of complex dimension $m$ that injects into $\mathrm{Aut}(N^\epsilon_\Lambda,J)$. The result then follows from Proposition 3 of \cite{istrati2019existence}: this torus has real dimension $2m\ge 4$ since $m\ge 2$. Moreover, it is stable under $J$. Thus, $N^\epsilon_\Lambda$ cannot admit an lck metric.
\end{proof}

We will need, in order to generalize Theorem \ref{CorollairecoolLVM}, the following definition.

\begin{defn}
Let $(\epsilon,\Lambda)$ and $(\epsilon,\Lambda')$ be two good systems with the same fundamental set $\epsilon$. We say that they are homotopic if there exists a continuous map
\begin{align*}
H\colon [0,1] &\rightarrow (\BC^m)^n\\
s &\mapsto H(s)
\end{align*}
satisfying the following conditions:
\begin{enumerate}
\item $H(0)=\Lambda$;
\item $H(1)=\Lambda'$ (more precisely, $H(1)$ is any permutation of $\Lambda'$);
\item for all $s\in [0,1]$, $(\epsilon, H(s))$ is a good system.
\end{enumerate}
\end{defn}
This definition naturally generalizes the one in \cite{meersseman1998procede} for good systems. Thus, lets $N^\epsilon_\Lambda$ and $N^\epsilon_{\Lambda'}$ two LVMB manifolds where $(\epsilon,\Lambda)$ and $(\epsilon,\Lambda')$ are homotopic, then $N^\epsilon_\Lambda$ and $N^\epsilon_{\Lambda'}$ are diffeomorphic (this is a direct consequence of an immediate adaptation of Proposition 1 from the fifth part of \cite{meersseman1998procede} to this definition).

\begin{thm}\label{thm2}
Let $N^\epsilon_\Lambda$ be an LVMB manifold. Then $N^\epsilon_\Lambda$ is lck with potential if and only if $m=1$ and $k=2$ (in other words, if and only if $N^\epsilon_\Lambda$ is biholomorphic to a diagonal Hopf manifold). 
\end{thm}
\begin{proof}
Suppose that $N^\epsilon_\Lambda$ admits a lck metric with potential. If $m=1$, then by Theorem \ref{thm1}, $N^\epsilon_\Lambda$ can only be a diagonal Hopf manifold. Now suppose that $m \ge 2$. In this case, $N^\epsilon_\Lambda$ does not satisfy $(K)$. Moreover, the lck metric with potential is a structure stable under deformation. Therefore, there exists a homotopic deformation of $(\epsilon,\Lambda)$ to $(\epsilon,\Lambda')$ such that $N^\epsilon_\Lambda$ and $N^\epsilon_{\Lambda'}$ are diffeomorphic, that $N^\epsilon_{\Lambda'}$ retains its lck with potential property, and that $N^\epsilon_{\Lambda'}$ satisfies condition $(K)$ since it is dense. Thus, Theorem \ref{thm1} forces the parameters of $N^\epsilon_{\Lambda'}$ to be $m=1$ and $k=2$. But since $N^\epsilon_\Lambda$ and $N^\epsilon_{\Lambda'}$ are diffeomorphic, the same holds for the parameters of $N^\epsilon_\Lambda$. It follows from Theorem \ref{caraHopf} that $N^\epsilon_\Lambda$ is biholomorphic to a diagonal Hopf manifold.
\end{proof}

\noindent We obtain the same conclusions as in \cite{faucard2024lvm} for LVMB manifolds (Table \ref{tab1}).
\begin{table}[h!]
\begin{tabular}{|c|c|c|c|}
\hline
Condition $(K)$ & $m$  & $k$ & Results on $N^\epsilon_\Lambda$ \\
\hline
& $1$  & $0$ or $1$ & $N^\epsilon_\Lambda$ not lck \\
\hline
& $1$  & $2$ & $N^\epsilon_\Lambda=$ Diagonal Hopf \\
\hline
& $\ge 2$ & any & $N^\epsilon_\Lambda$ not lck with potential, $N^\epsilon_\Lambda$ lck? \\
\hline
satisfied & $\ge 2$ & any & $N^\epsilon_\Lambda$ is not lck \\
\hline
\end{tabular}
\caption{Summary of results for an LVMB manifold $N^\epsilon_\Lambda$ with $m\ge 1, n>2m+1$ and $k$ indispensable coordinates.}
\label{tab1}
\end{table}

\section{Some LVMB manifolds are $(m-1)$-lck with potential}

We now generalize Theorem 6.2 from \cite{faucard2024lvm} by being more precise about the condition, which we will denote here as $(H)$.

\begin{defn}
Let $ N^\epsilon_\Lambda $ be an LVMB manifold with $ m \geq 2, n > 2m+1 $, and suppose that the coordinates $1,\ldots,m+1$ are indispensable. We say that $ N^\epsilon_\Lambda $ satisfies the condition $(H)$ if there exists a basis $ \mathcal{B} = (f_1, \ldots, f_m) $ of $ \BZ^m $ and $ j \in \{1, \ldots, m\} $ such that for all $ r \in \{m+2, \ldots, n\} $,
\[
|\Gamma(f_j)_r| < 1.
\]
\end{defn}

\begin{rem}
If such a $ j $ exists, then the action of $ \BZ $ generated by $ \alpha_j^\mathcal{B} $ is an action by holomorphic contractions. In the case where $ m = 1 $, $ n > 3 $, and the only indispensable coordinates are 1 and 2, the condition $(H)$ is always satisfied: this is what we showed in $(3) \Rightarrow (1)$ of Theorem \ref{caraHopf}.
\end{rem}

\begin{thm}\label{pstructure}
Let $ N^\epsilon_\Lambda $ be of type $ m \geq 2 $, $ n > 2m+1 $, and with $ k \geq m+1 $ essential points, the first $ m+1 $ coordinates being indispensable. Suppose that $ N^\epsilon_\Lambda $ satisfies condition $(H)$ for some basis $ \mathcal{B} = (f_1, \ldots, f_m) $ of $ \BZ^m $ and some $ j \in \{1, \ldots, m\} $. Then $ C/\langle \alpha_j^\mathcal{B} \rangle $ is a non-compact, non-Kählerian lck with potential $ \BZ^{m-1} $-covering of $ N^\epsilon_\Lambda $.
  
Moreover, since condition $(H)$ is open, this result is stable under homotopic deformation of the good system $ (\epsilon, \Lambda) $.
\end{thm}

\begin{proof}
The proof is exactly the same as that of Theorem 6.2 or Remark 6.4 from \cite{faucard2024lvm}. Remark 6.5 from \cite{faucard2024lvm} shows that the resulting covering is non-Kählerian.
\end{proof}
 
We immediately obtain the following corollary.

\begin{cor}\label{corsousfamille}
Keep the notation from Theorem \ref{pstructure}. Suppose that condition $(H)$ is satisfied for some basis $\mathcal{B}=(f_1,\ldots,m\}$ for a subfamily of indices $j_1<j_2<\cdots <j_l$ with $1\le l\le m$ and $j_i \in \{1,\ldots,m\}$ for all $i\in \{1,\ldots,l\}$. Then $C/\langle \alpha_{j_1}^\mathcal{B},\ldots,\alpha_{j_l}^\mathcal{B}\rangle$ is a non-compact, non-Kähler, lck with potential, $\BZ^{m-l}$-covering of $N^\epsilon_\Lambda$ (and this result is still stable under homotopic deformation).
\end{cor}

Before presenting some examples, let's recall the terminology introduced in \cite{faucard2024lvm}, which is more suitable for studying the existence of such coverings.

\begin{defn}
Let $M$ be a compact complex manifold. Let
\begin{align*}
p=\min\{&q\in \BN, \text{ there exists a connected, Galois, and lck $\BZ^q$-covering of } M\\
&\text{(respectively lck with potential) of } M\}.
\end{align*}
We say that $M$ is $p$-lck (respectively $p$-lck with potential) or that $M$ admits a $p$-lck structure (respectively lck with potential). If no such covering exists for $M$, then $p=\min \emptyset = -\infty$.
\end{defn}

\begin{rem}\label{conj}
The conjecture formulated in \cite{faucard2022lckLVM} for LVM manifolds can reasonably also be stated in the case of LVMB manifolds: under the conditions of Theorem \ref{pstructure}, $N^\epsilon_\Lambda$ is $(m-1)$-lck with potential. In other words, if this conjecture is true, then in Corollary \ref{corsousfamille}, the integer $l$ is always equal to $1$. We will show that this conjecture is false in the next section.
\end{rem}

\section{Examples and conjecture}\label{exemples}

In this section, I present examples of good systems that satisfy the conditions of Theorems \ref{thm1} and \ref{pstructure}. Some of these will particularly help to refute the conjecture stated in Remark \ref{conj}. I will then formulate two new conjectures at the end of these examples, which are related to Theorem \ref{pstructure}.

\vspace{0.5cm}

In his paper \cite{bosio2001varietes}, Bosio shows that there exist certain LVMB manifolds whose good system is not that of an LVM manifold. This study is conducted in the case $m=1$. This increases, compared to \cite{faucard2024lvm}, the class of manifolds that satisfy the conditions of Theorem \ref{thm1} for $m=1$, and are therefore non-lck. However, Bosio leaves the following problem open:
\begin{center}
Is there a LVMB manifold whose good system is not the good system of an LVM manifold, and which is not biholomorphic to any LVM manifold?
\end{center}
The answer to this question is affirmative. S. Cupit-Foutout and D. Zafran proved this in \cite{cupit2007non}. The authors prove the following theorem (Theorem 3.6 in the paper \cite{cupit2007non}):
\begin{thm}
Let $N_\Lambda^\epsilon$ be an LVMB manifold with good system $(\epsilon, \Lambda)$ such that $\Lambda_i \neq \Lambda_j$ for all $i \neq j$. Suppose additionally that $N_\Lambda^\epsilon$ satisfies condition $(K)$, that $(\epsilon, \Lambda)$ is not the good system of an LVM manifold, and that
\begin{align*}
\bigcap_{\sigma \in \epsilon} \sigma = \emptyset.
\end{align*}
Then $N_\Lambda^\epsilon$ is not biholomorphic to any LVM manifold.
\end{thm}
Thus, Example 1.2 in \cite{cupit2007non} provides an LVMB manifold that is not biholomorphic to any LVM manifold. This example is a case where $m=1$, so Theorem \ref{thm1} applies, and it provides an even more general application: we have an example of a non-lck LVMB manifold that is not biholomorphic to any LVM manifold. Another natural question arises:
\begin{center}
Is there an LVMB manifold whose good system is not that of an LVM manifold, and which is not homeomorphic to any LVM manifold?
\end{center}
This remains an open problem, and some exploration paths are presented by J. Tambour in \cite{tambour2010complexes}.

\vspace{0.5cm}

Now, concerning the examples of applications of Theorem \ref{pstructure}, it is necessary to assume that $m \geq 2$, which introduces a much more computational aspect than a geometric one. Computational tools have supported me in the search for such examples (see Appendix \ref{Annexe1}). Examples of good systems satisfying condition $(K)$ are not lacking, as it is sufficient for $\Lambda \in ((\BQ[i])^m)^n$. However, I have not found an example of such a good system where the coordinates $1,\ldots,m+1$ are indispensable and where I can prove that it does not correspond to the good system of an LVM manifold. Furthermore, note that if the cardinality of the fundamental set $\epsilon$ is $2$, then $(\epsilon, \Lambda)$ is the good system of an LVM manifold, according to the imbrication condition and Theorem \ref{caractLVM}. It must then be that $|\epsilon| \geq 3$ to have a chance to exhibit a LVMB manifold that is not an LVM manifold. This cardinality condition is equivalent to the inequality $n > 2m + 2$. Indeed, if $1, \ldots, m+1$ are indispensable, then
\begin{align*}
\epsilon = \{(1, \ldots, m+1, \sigma'), \sigma' \in \epsilon'\},
\end{align*}
where the elements of $\epsilon'$ have cardinality $m$. Moreover, the cardinality $e$ of $\epsilon'$ (and thus of $\epsilon$) is bounded above by the number of $m$-element subsets among $n-m-1$ elements: $\begin{pmatrix} n-m-1 \\ m \end{pmatrix}$. The $(PEUR)$ further lowers this bound. Indeed, if $\sigma' \in \epsilon$, then for all $k \in \{m+2, \ldots, n\} \setminus \sigma'$, there must exist a unique $k' \in \sigma'$ such that $\sigma' \setminus \{k'\} \cup \{k\} \in \epsilon'$. Fix such a $k$, and the set
\begin{align*}
\{k' \in \sigma', \sigma' \setminus \{k'\} \cup \{k\} \in \epsilon'\}
\end{align*}
has cardinality $1$. In other words, this shows that
\begin{align*}
e \leq \begin{pmatrix} n-m-1 \\ m \end{pmatrix} - (m-1).
\end{align*}
Therefore, if $n = 2m+2$, we have
\begin{align*}
\begin{pmatrix} n-m-1 \\ m \end{pmatrix} = \begin{pmatrix} m+1 \\ m \end{pmatrix} = m+1,
\end{align*}
and then $e \leq m+1 - (m-1) = 2$. Since $e > 1$ (otherwise the $(PEUR)$ condition is not satisfied), we have $e = 2$. I present several examples starting from $m=2$ and $n=7$ to test the conjecture in Remark \ref{conj}.

\begin{ex}\label{m2n7}{$m=2, n=7$}\\
The set
\begin{align*}
\epsilon=\{(12346),(12356),(12347),(12357)\}
\end{align*}
is fundamental. Consider the configuration
\begin{align*}
\Lambda=\begin{pmatrix}
0 & -1 & 0 & 1 & 1 & -1-i & -1-i\\
-i & 0 & -1 & i & i & 1 & 1 \end{pmatrix}.
\end{align*}
It is easy to verify (computationally, see Appendix \ref{Annexe1}) that $(\epsilon, \Lambda)$ is a good system.

\vspace{0.5cm}

Let us check that condition $(H)$ is satisfied. Let
\begin{align*}
\mathcal{B}=(f_1,f_2)=\Big(\begin{pmatrix}0\\ -2\end{pmatrix},\begin{pmatrix} 1\\ 0\end{pmatrix}\Big),
\end{align*}
which is a basis of $\BZ^2$. It remains to specify $\Gamma(f_1)$ and $\Gamma(f_2)$. The matrix $A$ is written as:
\begin{align*}
A=\begin{pmatrix} \Lambda_2-\Lambda_1\\ \Lambda_3-\Lambda_1\end{pmatrix}=\begin{pmatrix} -1 & i \\ 0 & -1+i\end{pmatrix}
\end{align*}
We have $\det A = 1-i \neq 0$, so $A$ is invertible with inverse:
\begin{align*}
A^{-1}=\frac{1}{1-i}\begin{pmatrix} -1+i & -i \\ 0 & -1\end{pmatrix}=\frac{1}{2}\begin{pmatrix} -2 & 1-i \\ 0 & -1-i\end{pmatrix}
\end{align*}
Thus,
\begin{align*}
\gamma(f_1)=2i\pi A^{-1}f_1=-2\pi\begin{pmatrix} 1+i \\ 1-i\end{pmatrix}
\end{align*}
while $\gamma(f_2)=\pi\begin{pmatrix} -2i\\ 0\end{pmatrix}$. We have
\begin{align*}
\Lambda_4-\Lambda_1=\Lambda_5-\Lambda_1=\Lambda_4-\Lambda_1=\begin{pmatrix}1\\ 2i\end{pmatrix}
\end{align*}
since $\Lambda_4=\Lambda_5$. Similarly,
\begin{align*}
\Lambda_6-\Lambda_1=\Lambda_7-\Lambda_1=\begin{pmatrix} -1-i\\ 1+i \end{pmatrix}
\end{align*}
And then
\begin{align*}
|\Gamma(f_1)_4|&=|\Gamma(f_1)_5|=|e^{-6\pi-6i\pi}|=e^{-6\pi}<1,\\
|\Gamma(f_1)_6|&=|\Gamma(f_1)_7|=|e^{4i\pi-4\pi)}|=e^{-4\pi}<1.
\end{align*}
Therefore, condition $(H)$ is satisfied here: $C/\langle \alpha_1^\mathcal{B}\rangle$ is a non-compact, non-Kähler, lck with potential $\BZ$-covering of $N_\Lambda^\epsilon$. Let us add two clarifications: on one hand, $C=(\BC^2\backslash\{0\})^2$, and on the other hand,
\begin{align*}
|\Gamma(f_2)_4|=|e^{-2i\pi}|=1.
\end{align*}
Thus, the hypothesis $(H)$ is satisfied only for $f_1$, but not for $f_2$. Consequently, Theorem \ref{pstructure} only tells us that $N_\Lambda^\epsilon$ is $p$-lck with potential, with $p \in \{0,1\}$. Furthermore, since the coordinates of the $\Lambda_j$ are rational, $N_\Lambda^\epsilon$ satisfies condition $(K)$ and is therefore non-lck. This implies that $p=1$. 

\vspace{0.5cm}

There remains one question: is it an LVM manifold? The answer to this question is yes. Indeed, if 
\begin{align*}
\theta=\{z\in \BC^4, \forall P \subset \{1, \ldots, 7\}, |P|=4, z \notin \mathcal{H}(\Lambda_j)_{j\in P}\}
\end{align*}
and if $O$ denotes the interior of $\mathcal{H}(\Lambda_j)_{j\in \sigma}$ (which is the same set for all $\sigma \in \epsilon$), then
\begin{align*}
\epsilon=\{\sigma\in \mathcal{P}(\{1, \ldots, 7\}), |\sigma|=5, O \subset \mathcal{H}(\Lambda_j)_{j\in \sigma}\}
\end{align*}
This shows (Theorem \ref{caractLVM}) that $N_\Lambda^\epsilon$ is an LVM manifold. Furthermore, one can show by considering the projections that it is biholomorphic to the product of two diagonal Hopf manifolds (thus indeed non-lck), and it corresponds to Example 1 in \cite{faucard2024lvm}.
\end{ex}

\begin{ex}\label{m2n8}{$m=2, n=8$}\\
Consider the fundamental set
\begin{align*}
\epsilon=\{(12345),(12346),(12357), (12348),(12367),(12378)\},
\end{align*}
and the configuration
\begin{align*}
\Lambda=\begin{pmatrix}
-2i & 1+i & -2-i & -2+i & -1+i & -1+i & -1+i & i \\
-i & -2 & -1-2i & -2-2i & 1-2i & -2-i & -2-2i & 1+i
\end{pmatrix}
\end{align*}
Then, $(\epsilon, \Lambda)$ is a good system that satisfies all the conditions of Theorem \ref{thm1}, so $N_\Lambda^\epsilon$ is non-lck. Moreover, condition $(H)$ is also satisfied for the basis
\begin{align*}
\mathcal{B}=(f_1,f_2)=\begin{pmatrix}
-1 & 0 \\ 2 & -1\end{pmatrix}
\end{align*}
and for $j=2$ (but not for $j=1$). In other words, according to Theorem \ref{pstructure}, $N_\Lambda^\epsilon$ is $1$-lck and $1$-lck with potential. The non-compact, non-Kähler, lck with potential covering is $C/\langle \alpha_2^\mathcal{B}\rangle$.
\end{ex}

\begin{ex}\label{m3n9}{$m=3, n=9$}\\
Consider the fundamental set
\begin{align*}
\epsilon=\{(1234567),(1234569),(1234589),(1234678),(1234789)\},
\end{align*}
and the configuration 
\begin{align*}
\Lambda=(\Lambda_1,\ldots,\Lambda_9)
\end{align*}
where
\begin{align*}
(\Lambda_1,\Lambda_2,\Lambda_3,\Lambda_4)&=\begin{pmatrix}
1+i & 1 & i & 1-2i\\
1-2i & 1+i & 1+i & 1-2i \\
0 & 1-2i & -1+i & -1-i
\end{pmatrix}\\
(\Lambda_5,\Lambda_6,\Lambda_7,\Lambda_8,\Lambda_9)&=\begin{pmatrix}
-1+i & 0 & -2i &-2-i & -2i \\
-2+i & -2i & -1 & -1-i & -1-i \\
1 & 1 & -1-2i & -2-2i & i
\end{pmatrix}.
\end{align*}
Then, $(\epsilon, \Lambda)$ is a good system satisfying condition $(K)$, so the associated manifold $N_\Lambda^\epsilon$ is non-lck. Moreover, condition $(H)$ is satisfied for the basis 
\begin{align*}
\mathcal{B}=(f_1,f_2,f_3)=\begin{pmatrix} -2 & 0 & 1 \\ 0 & 0 & -2 \\ 1 & -1 & -2\end{pmatrix}
\end{align*}
with $j=2$. It follows from Theorem \ref{pstructure} that $N_\Lambda^\epsilon$ is $p$-lck and $p$-lck with potential with $p\in \{1,2\}$. Indeed, since $N_\Lambda^\epsilon$ is non-lck, we have $p\ne 0$, and Theorem \ref{pstructure} implies that there exists an lck with potential non-compact, non-Kähler $\BZ^2$-covering of $N$. However, the conjecture in Remark \ref{conj} is refuted. Indeed, we have
\begin{align*}
&|\Gamma(f_2)_r|<1,\; \forall r\in \{5,6,7,8,9\}\\
&|\Gamma(f_3)_r|<1,\; \forall r\in \{5,6,7,8,9\}.
\end{align*}
But these modules are all strictly greater than $1$ for $f_1$. Thus, according to Corollary \ref{corsousfamille}, there exists an lck with potential non-compact, non-Kähler $\BZ$-covering of $N_\Lambda^\epsilon$: this is $C/\langle \alpha_2^\mathcal{B}, \alpha_3^\mathcal{B}\rangle$, so $p=1$. It is easy to generate other bases than $\mathcal{B}$ for which the same conclusion can be drawn (see Appendix \ref{Annexe1}).
\end{ex}

\begin{ex}\label{m3n10}{$m=3, n=10$}\\
Consider the fundamental set
\begin{align*}
\epsilon=\{&(1234567),(1234678),(1234579),(1234789),\\
&(123456\;10),(123468\; 10),(123459\; 10),(123489\; 10)\}
\end{align*}
and the configuration
\begin{align*}
\Lambda=(\Lambda_1,\ldots,\Lambda_{10}),
\end{align*}
where
\begin{align*}
(\Lambda_1,\Lambda_2,\Lambda_3,\Lambda_4,\Lambda_5)&=\begin{pmatrix}-1+2i & -2-i & 1     & 2     & -2   \\2+2i & 2    & 2-2i & 1     & -1-2i \\1&1-2i & 2-i   & 2-2i & 1    \end{pmatrix}\\
(\Lambda_6,\Lambda_7,\Lambda_8,\Lambda_9,\Lambda_{10})&=\begin{pmatrix}
 -2-2i  & -i & 0 & -1-2i & 2-2i \\
   -2      & -2 & -2-2i & -2+i & -2\\
          1-i     & -1+i & -1+2i & 2 & -2-2i\end{pmatrix}.
\end{align*}

Condition $(H)$ is satisfied for the $\BZ$-basis
\begin{align*}
\mathcal{B}=(f_1,f_2,f_3)=
\begin{pmatrix}
0 & 1 &0 \\
-2 & -1 & 0 \\
-2 & -2 & 1\end{pmatrix}
\end{align*}
for $j=1$ and $j=2$. Therefore, we can draw the same conclusion for this example as in Example \ref{m3n9}.
\end{ex}

\begin{ex}\label{m4n11}{$m=4, n=11$}\\
Consider the fundamental set 
\begin{align*}
\epsilon=\{(123456789),(12345789\; 10),(12345689\; 11),(1234589\; 10\; 11)\},
\end{align*}
and the configuration
\begin{align*}
\Lambda=(\Lambda_1,\ldots,\Lambda_{11}),
\end{align*}
where
\begin{align*}
(\Lambda_1,\Lambda_2,\Lambda_3)&=\begin{pmatrix} -1-i & -1 & -2 \\ -1-2i & 1-i & -i \\1 & i & -i \\ -2-i & 1-i & -i\end{pmatrix}\\
(\Lambda_4,\Lambda_5,\Lambda_6)&=\begin{pmatrix} -2-2i & -2+i & -i \\ -2-2i& -2 & i \\ 0& -1+i & -2i\\ 0 & -1-i & -2i\end{pmatrix}\\
(\Lambda_7,\Lambda_8,\Lambda_9)&=\begin{pmatrix} -1-i & -2-2i & -2+i \\ -2i & 1+i & -1-i \\ -1-i & 1-i & -2 \\ -2i & -2-2i & -1-2i\end{pmatrix}\\
(\Lambda_{10},\Lambda_{11})&=\begin{pmatrix} 1 & 1 \\ 1-2i & -2i \\ -1 & 1-i \\ -i & -2i\end{pmatrix}.
\end{align*}
Then, $(\epsilon,\Lambda)$ is a good system satisfying condition $(K)$ and condition $(H)$ for the basis
\begin{align*}
\mathcal{B}=(f_1,f_2,f_3,f_4)=\begin{pmatrix}-2 & -2 & 1 & 0 \\ -1 & 1 & 1 & 0 \\ 0 & 0& 0 & -1\\ -2 & 1 & 1 & 1\end{pmatrix}
\end{align*}
for $j=2,3$ and $4$. In other words, $N_\Lambda^\epsilon$ is $1$-lck with potential.
\end{ex}

We can summarize these examples in the following table (\ref{tab2}):
\begin{table}[h!]
\begin{center}
\begin{tabular}{|ccccc|}
\hline
Example & $m$ & $n$ & Biholomorphic to a LVM & $1$-LCK with potential \\
\hline
1.2 of \cite{cupit2007non} & $1$ & $5$ & no &  ? (non LCK)\\
\hline
\ref{m2n7} & $2$ & $7$ & yes & yes\\
\hline
\ref{m2n8} & $2$ & $8$ & ? & yes\\
\hline
 \ref{m3n9}& $3$ & $9$ & ? & yes\\
\hline
\ref{m3n10} & $3$ & $10$ & ? & yes\\
\hline
 \ref{m4n11}& $4$ & $11$ & ? & yes\\
\hline
\end{tabular}
\end{center}
\caption{Summary of the presented examples}
\label{tab2}
\end{table}

Based on these examples, it is natural to formulate the conjecture \ref{conjecture0}.

\newpage
\begin{appendix}

\section{Python scripts used for the search of examples}
\label{Annexe1}
The necessary packages:
\begin{lstlisting}[language=python]
import numpy as np
from scipy.optimize import linprog
import sympy as sp
from sympy.solvers import solve
\end{lstlisting}

The fundamental set $\epsilon$ and the configuration $\Lambda$ should be entered as global variables in the form of a list of lists. Concerning the configuration $\Lambda$, it is defined as an element of $(\BR^{2m})^n$, i.e., as a list of $n$ lists, each containing $2m$ elements. For example, if the first vector of the configuration $\Lambda$ is $\begin{pmatrix} -1+2i \\ -1+i\\ 1\end{pmatrix}$, then \verb?Lambda[0]? will be the list $[-1,2,-1,1,1,0]$. Let us denote by$N_\Lambda^\epsilon$the manifold associated with $(\epsilon, \Lambda)$ when this system is good.

 \subsection{Studyability of $(\epsilon,\Lambda)$}\label{AnnexeEtudiable}

 Let $\Delta_j = \Lambda_j - \Lambda_1$ for all $j \in \{2, \ldots, 10\}$. The function \verb?Delta(L)? returns the configuration where this subtraction by $\Lambda_1$ is performed.

  \begin{lstlisting}[language=python]
  def Delta(L):
    DeltaRes=[]
    for i in range(len(L)-1):
        DeltaRes.append([b - a 
        for a, b in zip(L[0],L[i+1])])
    return DeltaRes
  \end{lstlisting}

  The function \verb?is_free(L, sigma)? returns $1$ if $(\Delta_j)_{j \in \sigma \backslash \{1\}}$ is linearly independent and $0$ otherwise. The input \verb?sigma? corresponds to a certain fundamental set excluding $1$: $\sigma \backslash \{1\}$. Note that imposing an error of $0.001$ on the determinant accounts for the numerical nature of \verb?numpy? operations.

   \begin{lstlisting}[language=python] 
  def is_free(L,sigma):
    Delta1=Delta(L)
    Delta_sigma=[]
    for i in sigma:
        Delta_sigma.append(Delta1[i-2]) 
    Delta_sigma=np.array(Delta_sigma)
    a=np.linalg.det(Delta_sigma)
    if np.abs(a)>0.001: 
        return 1
    else:
        return 0
    \end{lstlisting}
    
It is then necessary to test all the fundamental subsets. \\
The function \verb?is_studyable(epsilon, L)? returns $1$ if the input configuration is studyable.

\begin{lstlisting}[language=python] 
def is_studyable(epsilon,L):
    N=len(epsilon)
    for i in range(N):
        sigma=epsilon[i]
        if is_free(L,sigma)==0:
            return 0
    return 1
\end{lstlisting}

 \subsection{Imbrication condition}\label{AnnexeImbrication}

Let $A$ be a subset of $\epsilon$. Let $H_A$ be the global intersection of the interiors of the convex hulls of the sub-configurations associated with the fundamental parts of $A$. Let $a = |A|$. Then, if we define $a \times (2m+1)$ strictly positive variables, say $a_j^\sigma > 0$ for $j = 1, \ldots, 2m+1$ and $\sigma \in A$, the set $H_A$ is non-empty if the system
\begin{align*}
\begin{cases}
\sum_{i \in \sigma} a^\sigma_i &= 1 \quad \forall \sigma \in A, \\
\sum_{i \in \sigma} a^\sigma_i \Lambda_i &= \sum_{i \in \tau} a^\tau_i \Lambda_i \quad \forall (\sigma, \tau) \in A^2
\end{cases}
\end{align*}
has at least one solution. The function \verb?is_empty(A,L)? returns $1$ if the system associated with $A$ has no solution, after having declared the $a \times (2m+1)$ positive variables of the system via \verb?sympy?. The imbrication condition will be satisfied if $H_A$ is non-empty for every subset $A$ of $\epsilon$ of cardinality $2$.

 \begin{lstlisting}[language=python]  
def symbolic_variables(m,a): 
    V = []
    for i in range((2*m+1)*a):
        name = 'x' + str(i)
        name=sp.symbols(name,positive=True)
        V.append(name)
    return V

  
def is_empty(A,L):
    a=len(A)
    n=len(L)
    m=int(len(L[0])/2)
    V=symbolic_variables(m,a)
    Eq=[]
    for i in range(a):
        eq=-1
        for j in range(2*m+1):
            eq+=V[(2*m+1)*i+j]
        Eq.append(eq)
    for k in range(a-1):
        for i in range(2*m):
            eq=0
            for l in range(2*m+1):
                eq+=L[A[k][l]-1][i]*V[l+(2*m+1)*k]
            for j in range(2*m+1):
                eq=eq-L[A[-1][j]-1][i]*V[j+(2*m+1)*(a-1)]
            Eq.append(eq)
    Solutions=solve(Eq,V)
    if []==Solutions:
        return 1
    return 0


def Imbrication_condition(epsilon,L):
    n_eps=len(epsilon)
    i=0
    while i<n_eps:
        j=i+1
        while j<n_eps:
            if is_empty([epsilon[i],epsilon[j]],L))
                return 0
            j+=1
        i+=1
    return 1

  \end{lstlisting}

 \subsection{Is $(\epsilon,\Lambda)$ the good system of an LVM manifold ?}\label{AnnexeLVM}
 
 If there exists a subset $A$ of $\epsilon$ with cardinality strictly greater than $2$ such that $H_A$ is empty, then theorem \ref{caractLVM} implies that $N_\Lambda^\epsilon$ is not an LVM manifold. I then attempted, for each good system $(\epsilon,\Lambda)$, to find such an $A$ by testing all subsets of $\epsilon$ each time, but without success.

\begin{lstlisting}[language=python]
def partlist(seq):
    p = []
    imax = 2**len(seq)-1
    for i in range(imax):
        s = []
        j, jmax = 0, len(seq)-1
        while j <= jmax:
            if (i>>j)&1 == 1:
                s.append(seq[j])
            j += 1
        p.append(s)
        i += 1 
    return p

def is_LVM(epsilon,L):
    l=len(epsilon)
    n=len(L)
    m=int(len(L[0])/2)
    SS_ENS=partlist(epsilon)
    nb_part=len(SS_ENS)
    for i in range(nb_part):
        A=SS_ENS[i]
        if len(A)>2 and is_empty(A,L):
            return 0
    return 1
    \end{lstlisting}

\subsection{Condition $(H)$}

This test requires the following three functions. The first, \verb?Lambda_C(L)? for a configuration \verb?L? as an element of $(\BR^{2m})^n$, returns the configuration \verb?L? but as an element of $(\BC^m)^n$. The second, \verb?Gamma(k,r,L)? takes as input an element $k \in \BZ^m$, an integer $r \in \{m+2, \ldots, n\}$, as well as a configuration \verb?L? as an element of $(\BR^{2m})^n$ and returns $\Gamma(k)_r$. Finally, the third function returns $1$ if condition $(H)$ is satisfied for a certain base \verb?B? for the good system \verb?epsilon,L?.

\begin{lstlisting}[language=python]
def Lambda_C(L):
    n_L=len(L)
    m=int(len(L[0])/2)
    Res=[]
    for i in range(n_L):
        Res.append([])
        for j in range(m):
            Res[i].append(L[i][2*j]+L[i][2*j+1]*1j)
    return Res

def Gamma(k,r,L):
    Delta1=Delta(L,1)
    Delta1_C=Lambda_C(Delta1)
    m=int(len(L[0])/2)
    A=[]
    for i in range(m):
        A.append(Delta1_C[i])
    A=np.array(A)
    InvA=np.linalg.inv(A)
    Right=2*1j*np.pi*np.dot(InvA,k)
    Left=Delta1_C[r-2]
    return np.exp(np.dot(Left,Right))

def Cond_H(B,epsilon,L):
    n=len(L)
    m=int(len(L[0])/2)
    Res=[]
    for i in range(m):
        Res.append([])
        for r in range(m+2,n+1):
            Res[i].append(np.abs(Gamma(B[i],r,L)))
    j=[]
    liste_max=[]
    for i in range(m):
        m=max(Res[i])
        if m<1:
            j.append(i+1)
    if j==[]:
        return 0
    else:
        print(j)
        return 1
\end{lstlisting}

\subsection{Randomly generate configurations $\Lambda$}
For the search of examples, the random generation of configurations proves to be useful. The function \verb?random_Conf_R(m,n)? randomly generates an element of $(\BR^{2m})^n$. The function \verb?generate_good_system(epsilon)? provides a configuration \verb?Lambda? such that \verb?(epsilon,Lambda)? is a good system.

\begin{lstlisting}[language=python]
def random_Conf_R(m,n):  
    L=[]
    for j in range(n):
        L.append([])
        for i in range(2*m):
            x=np.random.choice(np.arange(-2,2,1))
            L[j].append(x)
    return L
    
def generate_good_system(epsilon): 
    M=[]
    for i in range(len(epsilon)):
        M.append(max(epsilon[i]))
    n=max(M)
    m=int((len(epsilon[0])-1)/2)
    L=random_Conf_R(m,n)
    cpt=1
    while not(Imbrication_condition(epsilon,L) and 
    is_studyable(epsilon,L)):
        L=random_Conf_R(m,n)
        cpt+=1  
        if cpt>10000:
            print("too many trials")
            return None
    print("number of trials=",cpt)
    return L
 \end{lstlisting}
 
 Once a good system is available, one can generate a $\BZ^m$-basis such that the condition $(H)$ is satisfied. It turns out that this function has always returned such a basis, and this is why I formulated conjecture \ref{conjecture0}.

\begin{lstlisting}[language=python]
def base(m):
    B=[]
    for j in range(m):
        B.append([])
        for i in range(m):
            x=np.random.choice(np.arange(-2,2,1))
            B[j].append(x)
    B=np.array(B)
    while np.linalg.det(B)<0.01:
        B=[]
        for j in range(m):
            B.append([])
            for i in range(m):
                x=np.random.choice(np.arange(-2,2,1))
                B[j].append(x)
        B=np.array(B)
    return B
    
    
def generate_base_H(epsilon,L):
    m=int(len(L[0])/2)
    B=base(m)
    print(B)
    cpt=1
    while not(Cond_H(B,epsilon,L)):
        B=base(m)
        cpt+=1
        if cpt>=10000:
            print("too many trials=",cpt)
            return [0]
    print("numbre of trials=",cpt)
    return B
 \end{lstlisting}
\end{appendix}

\bibliographystyle{spmpsci}
\bibliography{myBibLib} 

\end{document}